\documentclass[a4paper,11pt,twoside,reqno]{amsart}

\usepackage[utf8]{inputenc}
\usepackage[plainpages=false,pdfpagelabels=true]{hyperref}
\usepackage{amssymb,amsthm}
\usepackage[margin=1in]{geometry}
\usepackage{slashed}

\newtheorem{Satz}{Theorem}[section]

\newtheorem{Lem}[Satz]{Lemma}

\theoremstyle{definition}
\newtheorem{Dfn}[Satz]{Definition}
\newtheorem{Bem}[Satz]{Remark}

\newcommand{\tr}{\operatorname{Tr}}

\newcommand{\vol}{{\operatorname{vol}}}
\newcommand{\dv}{\text{ }dv}
\allowdisplaybreaks[1]

\renewcommand{\epsilon}{\varepsilon}

\newcommand{\R}{\ensuremath{\mathbb{R}}}
\newcommand{\N}{\ensuremath{\mathbb{N}}}

\newcommand{\s}{\ensuremath{\mathbb{S}}}
\numberwithin{equation}{section}

\title{More weakly biharmonic maps from the ball to the sphere}
\author{Volker Branding}
\date{\today}
\address{University of Vienna, Faculty of Mathematics\\
Oskar-Morgenstern-Platz 1, 1090 Vienna, Austria\\}
\email{volker.branding@univie.ac.at}
\thanks{The author gratefully acknowledges the support of the Austrian Science Fund (FWF) through the project "Geometric Analysis of Biwave Maps" (DOI: 10.55776/P34853)
}
\begin{document}

\begin{abstract}
In this note we prove the existence of two proper biharmonic maps between
the Euclidean ball of dimension bigger than four and Euclidean spheres of appropriate dimensions.
We will also show that, in low dimensions, both maps are unstable critical points of the bienergy.
\end{abstract} 

\maketitle

\section{Introduction and results}
In the geometric calculus of variations one is interested in the construction
of non-trivial maps between two Riemannian manifolds \((M,g)\) and \((N,h)\) which are critical points of a given functional. In this regard, much attention has been paid to the case of the classic \emph{energy of a map} given by
\begin{align}
\label{energy}
E(\phi):=\frac{1}{2}\int_M|d\phi|^2\dv,
\end{align}
where \(\phi\colon M\to N\) is a map which we assume to be smooth for the moment.
The critical points of \eqref{energy} are characterized by the vanishing
of the so-called \emph{tension field}, i.e.
\begin{align*}
0=\tau(\phi):=\tr_g\bar\nabla d\phi,\qquad \tau(\phi)\in\Gamma(\phi^\ast TN).
\end{align*}
Here, \(\bar\nabla\) represents the connection on the pull-back bundle \(\phi^\ast TN\). The famous  result of Eells - Sampson \cite{MR164306} guarantees the existence of a harmonic map in each homotopy class of maps if both \(M,N\) are closed and if \(N\)
has non-positive sectional curvature. Note that these harmonic maps are stable in the sense that 
they reach the minimum of the energy in their homotopy class \cite{MR214004,MR375386}.
For more details and references on the stability of harmonic maps we refer to \cite{MR4477489}.

In the case of a spherical target the existence problem for harmonic maps is substantially more complicated.
In this setup it is convenient to consider \(\mathbb{S}^n\subset \R^{n+1}\) and the map \(u\colon M\to\mathbb{S}^n\subset\R^{n+1}\).
Then, a map to the sphere is harmonic if and only if it satisfies
\begin{align}
\label{eqn:harmonic-sphere}
\Delta u+|du|^2u=0,
\end{align}
where \(\Delta\) represents the Laplace-Beltrami operator on \(M\). In general, harmonic maps
to a spherical target are unstable, see \cite{MR4477489} and references therein for the precise details.
For an overview on harmonic maps we refer to the book \cite{MR2044031}.

In many nonlinear problems in geometric analysis it is favorable to not only consider the class of smooth maps but to also allow for maps of lower regularity and to consider a weak version of the problem at hand.
In order to approach the notion of weak solutions of the harmonic map equation
let us recall the definition of the Sobolev space for maps to the sphere
\begin{align*}
W^{p,q}(M,\s^n):=\{u\in W^{p,q}(M,\R^{n+1})~~\mid u(x)\in\s^n ~~ a.e.\}.    
\end{align*}
In the case that \(\partial M\neq 0\) we let \(u_0\in W^{p,q}(M,\s^n)\) and define
\begin{align*}
W^{p,q}_{u_0}(M,\s^n):=\{u\in W^{p,q}(M,\s^n) \mid \nabla^k(u-u_0)\big|_{\partial M}=0,
0\leq k\leq p-1\}.
\end{align*}
Here, the boundary condition is to be understood in the sense of traces.

For \(\Omega\subset M\) open we say that a map \(u\in W^{1,2}(\Omega,\s^n)\) is a \emph{weak harmonic map} if it solves \eqref{eqn:harmonic-sphere} in the sense of distributions.
A celebrated result of H\'{e}lein \cite{MR1078114} shows that such weak harmonic maps are actually smooth if the domain is two-dimensional 
which led to the general regularity theory
of Rivi\`ere \cite{MR2285745} for two-dimensional conformally invariant variational problems.
In \cite{MR2383929} Rivi\`ere and Struwe generalized this analysis to higher dimensions. 

A fourth order generalization of harmonic maps, that receives growing attention both in analysis and geometry, is given by the theory of \emph{biharmonic maps}.
Here, the starting point is the \emph{bienergy functional} which is given by
\begin{align}
\label{bienergy}
E_2(\phi):=\frac{1}{2}\int_M|\tau(\phi)|^2\dv.
\end{align}
The critical points of \eqref{bienergy} are characterized by the vanishing of the \emph{bitension field}
\begin{align}
\label{bitension}
0=\tau_2(\phi):=\bar\Delta\tau(\phi)+\tr_g R^N(\tau(\phi),d\phi(\cdot))d\phi(\cdot).    
\end{align}
In the above formula \(\bar\Delta\) represents the connection Laplacian on the pull-back bundle
\(\phi^\ast TN\).
A direct inspection of the biharmonic map equation reveals that every harmonic map automatically provides a solution. Hence, one is very much interested in finding the non-harmonic solutions of \eqref{bitension} which are called \emph{proper biharmonic}. However, in the case that \(M\) is closed and \(N\) has non-positive sectional curvature the maximum principle implies that every biharmonic map must be harmonic \cite{MR886529}. Further classification results for biharmonic maps can be found in \cite{MR4040175}.

Due to the reasons outlined above most research on biharmonic maps considers the case of a spherical target. In this setup the bienergy can be expressed as
\begin{align}
E_2(u)=\frac{1}{2}\int_M(|\Delta u|^2-|\nabla u|^4)\dv
\end{align}
and its critical points satisfy
\begin{align}
\label{eqn:biharmonic-sphere-extrinsic}
\Delta^2u
+2\operatorname{div}\big(|\nabla u|^2\nabla u\big)
+\big(|\Delta u|^2+\Delta|\nabla u|^2+2\langle\nabla u,\nabla\Delta u\rangle
+2|\nabla u|^4\big)u=0,
\end{align}
where again \(u\colon M\to\mathbb{S}^n\subset\R^{n+1}\).
Concerning the regularity of biharmonic maps to spheres we refer to
the articles \cite{MR1692148,MR2375314,MR2094320}, 
see also the more general result of Lamm and Rivi\`ere \cite{MR2383929} on the regularity
of fourth order problems in dimension four which was later extended by Struwe to higher dimensions
\cite{MR2413109}.

Unique continuation theorems for biharmonic maps, in particular for biharmonic maps to spheres, have been established in \cite{MR3990379}. Rotationally symmetric biharmonic maps between model spaces,
which include flat Euclidean space and the Euclidean sphere, have been investigated in \cite{MR3357596}.
For a detailed study of the stability of biharmonic maps to spheres we refer to \cite{MR2135286}, biharmonic homogeneous polynomial maps between spheres
have recently been constructed in \cite{MR4599666}.

The current status of research on biharmonic maps can be found in the recent book
\cite{MR4265170}, for more details on the general structure of higher order variational problems one can consult \cite{MR4106647} and references therein.

Before we state the main results of this article let us present a number of explicit solutions to the harmonic and biharmonic map equations \eqref{eqn:harmonic-sphere} and \eqref{eqn:biharmonic-sphere-extrinsic} for the case of a spherical target.

One prominent solution to \eqref{eqn:harmonic-sphere} is given by the 
equator map
\begin{align}
\nonumber
w\colon B^m&\to\mathbb{S}^m\subset \R^{m+1} \\
\label{sol:jk} x&\mapsto \frac{x}{r},
\end{align}
where \(B^m\) is the unit ball in \(\R^m\). 
Note that the equator map is a weak solution of the harmonic map equation \eqref{eqn:harmonic-sphere}.
We will use the notation
\begin{align*}
r:=\sqrt{x_1^2+\ldots+x_m^2},
\end{align*}
where \(x\in\R^m\) throughout the whole manuscript.

This particular solution of the harmonic map equation
was first mentioned by Hildebrandt in \cite[p. 13]{MR714341}.
Subsequently, Jäger and Kaul proved that this map is energy minimizing for \(m\geq 7\).
Later, Lin \cite{MR916327} realized that this map is energy minimizing
if interpreted as \(w\colon B^m\to\mathbb{S}^{m-1}\) for all \(m\geq 3\).

Now, we recall the following definition as it is central within this manuscript.
\begin{Dfn}
A map \(u\in W^{2,2}(B^m,\s^n)\) is called a weak biharmonic map if it solves \eqref{eqn:biharmonic-sphere-extrinsic} in the sense of distributions.
\end{Dfn}

Recently, Fardoun, Montaldo and Ratto showed in \cite[Theorem 1.1]{MR4076824} how to manufacture
a weak proper biharmonic map out of the equator map \eqref{sol:jk}. 
Their approach relies on "rotating" the harmonic map \eqref{sol:jk}
to a proper biharmonic one.
More precisely, they proved that the map
\begin{align}
\label{sol:fmr}
\nonumber
\tilde w\colon B^m&\to\mathbb{S}^m\subset \R^{m}\times\R \\
x&\mapsto \big(\sin a\frac{x}{r},\cos a\big)   
\end{align}
is a weak proper biharmonic map if and only if
\begin{enumerate}
    \item \(m=5\) and \(a=\frac{\pi}{3}\)
    \item \(m=6\) and \(a=\frac{1}{2}\operatorname{arccos}(-\frac{4}{5})\).
\end{enumerate}
They also proved that both of these proper biharmonic maps are unstable, see \cite[Theorem 1.2]{MR4076824}.

In another recent paper \cite{MR4371934} Misawa and Nakauchi 
found two generalizations of the equatorial harmonic map \eqref{sol:jk}.
In particular, they proved that the maps
\begin{align}
\nonumber
(u)_{ij}:=u_{ij}\colon B^m&\to \s^{m^2-1}\subset\R^{m^2},\\
\label{eqn:u} x&\mapsto\frac{1}{\sqrt{m(m-1)}}\big(-\delta_{ij}+m\frac{x_ix_j}{r^2}\big),\\
\nonumber
(v)_{ijk}:=v_{ijk}\colon B^m&\to \s^{m^3-1}\subset\R^{m^3},\\
\label{eqn:v} x&\mapsto\frac{1}{\sqrt{(m-1)(m+2)}}\big(
\delta_{ij}\frac{x_k}{r}+\delta_{jk}\frac{x_i}{r}
+\delta_{ik}\frac{x_j}{r}-(m+2)\frac{x_ix_jx_k}{r^3}\big)
\end{align}
are both harmonic, that is they solve the equation for harmonic maps to spheres \eqref{eqn:harmonic-sphere}. 
There is also a systematic, recursively defined, generalization of such harmonic maps to spheres, see \cite{MR4124860,MR4593065}
and their stability was recently investigated in \cite{MR4686854}.

We would like to point out that the initial motivation to study maps of
the form \eqref{eqn:u}, \eqref{eqn:v} in the 1980's by Giaquinta and Ne\v{c}as
was to find counterexamples to the regularity of weak solutions of elliptic systems
that were investigated at that time, see for example \cite{MR581329,MR1401417,MR566246}.

As both \(u_{ij}\) and \(v_{ijk}\) are harmonic maps to the sphere it seems very natural to ask if the same kind of "rotation" that turned the equatorial harmonic map \eqref{sol:jk} into a proper biharmonic map can also be applied
in the cases of the maps \(u_{ij}\) and \(v_{ijk}\). 

Our first result is to show that the answer is affirmative for \(u_{ij}\):

\begin{Satz}
\label{thm:biharmonic-u}
Let \(m>4\). The map 
\begin{align*}
(\tilde u)_{ij}:=\tilde u_{ij}\colon B^m\to\s^{m^2-1}\subset \R^{m^2-1}\times\R
\end{align*}
defined by
\begin{align}
\label{dfn:ualpha}
x\mapsto (\sin\alpha~u_{ij},\cos\alpha)    
\end{align}
is a proper weak biharmonic map if and only if 
\begin{align*}
\sin^2\alpha=1-\frac{2}{m},\qquad \alpha\in (0,\frac{\pi}{2}).
\end{align*}
\end{Satz}

A direct calculation shows that the proper biharmonic map provided by Theorem \ref{thm:biharmonic-u} has non-constant energy density 
\begin{align*}
\frac{|\nabla\tilde u_{ij}|^2}{2}=\frac{m-2}{r^2}    
\end{align*}
and bienergy
\begin{align*}
E_2(\tilde u_{ij})=4\frac{m-2}{m-4}\vol(\s^{m-1}).
\end{align*}
The last observation suggests that this map cannot be a stable critical point of the bienergy \eqref{eqn:biharmonic-sphere-extrinsic}.
Indeed, regarding the stability of the proper biharmonic map given by Theorem \ref{thm:biharmonic-u} we 
establish the following statement:
\begin{Satz}
\label{thm:biharmonic-u-stability}
The proper biharmonic map provided by Theorem \ref{thm:biharmonic-u}
is unstable if \(5\leq m\leq 12\).
\end{Satz}

In addition, we will also show how the map \(v_{ijk}\) defined in \eqref{eqn:v}
can be transformed into a proper biharmonic map.

\begin{Satz}
\label{thm:biharmonic-v}
Let \(m>4\). The map 
\begin{align}
(\tilde v)_{ijk}:=\tilde v_{ijk}\colon B^m\to\s^{m^3-1}\subset \R^{m^3-1}\times\R
\end{align}
defined by
\begin{align*}
\label{dfn:valpha}
x\mapsto (\sin\beta~v_{ijk},\cos\beta)    
\end{align*}
is proper biharmonic if and only if 
\begin{align}
\sin^2\beta=\frac{5}{6}\frac{m-1}{m+1},\qquad \beta\in (0,\frac{\pi}{2}).
\end{align}
\end{Satz}

Again, the biharmonic map provided by the previous Theorem has non-constant 
energy density, i.e.
\begin{align*}
\frac{|\nabla\tilde v_{ijk}|^2}{2}=\frac{5}{4}\frac{m-1}{r^2}.
\end{align*}
We will also investigate the stability of the proper biharmonic map constructed in the previous Theorem.
If we calculate its bienergy we find
\begin{align*}
E_2(\tilde v_{ijk})=\frac{5}{8}\frac{(m-1)(11m+1)}{m-4}\vol(\s^{m-1})
\end{align*}
again indicating that it cannot be a stable critical point of \eqref{bienergy}.
Indeed, we will prove the following
\begin{Satz}
\label{thm:biharmonic-v-stability}
The proper biharmonic map provided by Theorem \ref{thm:biharmonic-v}
is unstable if \(5\leq m\leq 18\).
\end{Satz}

Let us make the following remarks in order to put the results obtained in this manuscript into the bigger picture on biharmonic maps to spheres.

\begin{Bem}
\begin{enumerate}
    \item If we compare Theorems \ref{thm:biharmonic-u} and \ref{thm:biharmonic-v}
with the corresponding results on the proper biharmonic map that was obtained by rotating the equator map \eqref{sol:jk} then we find that we always need to require \(m>4\). However, rotating the equator map only gives rise to a proper biharmonic map if \(m=5,6\) while rotating the maps \(u_{ij}\)
and \(v_{ijk}\) provides proper biharmonic maps in all dimension bigger than four.

\item In \cite{MR4227593} Laurain and Lin proved a general existence result for biharmonic maps from the four-dimensional unit ball to the sphere using the heat flow method. While their general result applies to biharmonic maps in the critical dimension
the biharmonic maps obtained in Theorems \ref{thm:biharmonic-u} and Theorems \ref{thm:biharmonic-v} are defined in the supercritical case.

\item In order for the biharmonic maps \(\tilde u_{ij},\tilde v_{ijk}\) provided by Theorems \ref{thm:biharmonic-u} and \ref{thm:biharmonic-v} to be in the Sobolev space
\(W^{2,2}(B^m,\s^n)\) it is necessary that 
\begin{align*}
\int_0^1r^{m-5}dr<\infty.    
\end{align*}
It is easy to check that this condition is only satisfied if \(m>4\). Hence, Theorems \ref{thm:biharmonic-u} and \ref{thm:biharmonic-v} do not provide a solution
if \(m=3,4\) although the equations for \(\alpha,\beta\) can be solved in these cases.

\item One should expect that the biharmonic maps provided by Theorems \ref{thm:biharmonic-u} and \ref{thm:biharmonic-v} will be unstable in all dimensions. However, the method of proof that we employ in this paper only seems to be feasible in lower dimensions.

\item The biharmonic maps constructed in Theorem \ref{thm:biharmonic-u} (assuming \(m>2\)) and in Theorem \ref{thm:biharmonic-v} (assuming \(m>1\)) are also solutions to the biharmonic map equation if \(M=\R^m\setminus\{0\}\). In this case they are even smooth solutions of \eqref{eqn:biharmonic-sphere-extrinsic} as we have excluded the origin of \(\R^m\). However, we cannot study the stability of these biharmonic maps as the domain is non-compact and for this reason we have formulated our main results
for the case that the domain is \(B^m\).

\item In Theorem 18 of \cite{MR886529} it is shown that a biharmonic map from a closed Riemannian manifold to the sphere is unstable under the assumption 
that it has constant energy density and satisfies the so-called conservation law.
As the maps constructed in Theorems \ref{thm:biharmonic-u} and \ref{thm:biharmonic-v} do not have constant energy density and as the unit ball \(B^m\)
has a boundary the aforementioned result cannot be applied in our case.
\end{enumerate}

\end{Bem}

Throughout this article we will use the following sign conventions:
For the Riemannian curvature tensor field we use 
$$
R(X,Y)Z=[\nabla_X,\nabla_Y]Z-\nabla_{[X,Y]}Z,
$$ 
where \(X,Y,Z\) are vector fields.
For the rough Laplacian on the pull-back bundle $\phi^\ast TN$ we employ the analysts sign convention
$$
\bar\Delta = \tr_g(\bar\nabla\bar\nabla-\bar\nabla_\nabla).
$$
We will use Latin letters to represent indices on the domain \(M\).
Moreover, we will always employ the summation convention and tacitly sum over repeated indices.

\section{Proof of the main results}
In this section we provide the proofs of the main theorems.
First, we establish the following 
\begin{Lem}
\label{lem:der-u}
Let \(u_{ij}\colon B^m\to\s^{m^2-1}\) be the map defined in \eqref{eqn:u}. 
Then, the following identities hold
\begin{align}
\label{identities-u1}
|\nabla u_{ij}|^2=&\frac{2m}{r^2}, \\
\nonumber\Delta u_{ij}=&\frac{2m}{\sqrt{m(m-1)}}\big(\frac{\delta_{ij}}{r^2}-m\frac{x_ix_j}{r^4}\big) \\
\nonumber=&-\frac{2m}{r^2}u_{ij}
,\\
\nonumber\Delta^2 u_{ij}=&\frac{8m(m-2)}{\sqrt{m(m-1)}}\frac{1}{r^4}
\big(    
-\delta_{ij}+m\frac{x_ix_j}{r^2}
\big)
=\frac{8m(m-2)}{r^4}u_{ij},
\end{align}
where \(\delta_{ij}\) represents the Kronecker delta.
\end{Lem}
\begin{proof}
A direct calculation shows the following identity
\begin{align*}
\nabla_k u_{ij}=&\frac{m}{\sqrt{m(m-1)}}\big(
\frac{\delta_{ik}x_j}{r^2}+\frac{x_i\delta_{jk}}{r^2}
-2\frac{x_ix_jx_k}{r^4}\big),
\end{align*}    
taking the square then gives the first claim.

Regarding the second claim we differentiate once more and find
\begin{align*}
\nabla_l\nabla_k u_{ij}=&\frac{m}{\sqrt{m(m-1)}}\big(\frac{\delta_{ik}\delta_{lj}}{r^2}
-2\frac{\delta_{ik}x_jx_l}{r^4}
+\frac{\delta_{il}\delta_{jk}}{r^2}
-2\frac{\delta_{jk}x_ix_l}{r^4} \\
&-2\frac{\delta_{il}x_jx_k}{r^4}
-2\frac{\delta_{jl}x_ix_k}{r^4}
-2\frac{\delta_{lk}x_ix_j}{r^4}
+8\frac{x_ix_jx_kx_l}{r^6}\big)
\end{align*}
and taking the trace gives the second formula of the Lemma.

Finally, we calculate
\begin{align*}
\nabla_k\Delta u_{ij}=\frac{4mx_k}{r^4}u_{ij}-\frac{2m}{r^2}\nabla_ku_{ij}
\end{align*}
and differentiating once more we obtain
\begin{align*}
\nabla_l\nabla_k\Delta u_{ij}=&
\frac{4m\delta_{kl}}{r^4}u_{ij}
-\frac{16mx_kx_l}{r^6}u_{ij}
+\frac{4mx_k}{r^4}\nabla_lu_{ij}
+\frac{4mx_l}{r^4}\nabla_ku_{ij}
-\frac{2m}{r^2}\nabla_l\nabla_k u_{ij}.
\end{align*}
Now, a direct calculation shows that \(x\cdot \nabla u_{ij}=0\) and
the last formula follows by taking the trace again and using the second formula
of the Lemma.
\end{proof}

\begin{proof}[Proof of Theorem \ref{thm:biharmonic-u}]
First of all, we recall that we can rewrite the equation
for biharmonic maps to spheres in the following form
\begin{align}
\label{eq:el-short}
\Delta^2u+2\operatorname{div}\big(|\nabla u|^2\nabla u\big)    
-\big(\langle\Delta^2u,u\rangle-2|\nabla u|^4\big)u=0
\end{align}
which is equivalent to \eqref{eqn:biharmonic-sphere-extrinsic}.
This can easily be seen by applying \(\Delta^2\) to \(|u|^2=1\) giving
\begin{align*}
0=\Delta|\nabla u|^2+|\Delta u|^2+2\langle\nabla\Delta u,\nabla u\rangle
+\langle u,\Delta^2u\rangle.
\end{align*}
Now, we inspect the precise structure of the map \(\tilde u_{ij}\) defined
in \eqref{dfn:ualpha}.
Since the \(m^2\)-th component of the map \(\tilde u_{ij}\)
is constant, \eqref{eq:el-short} yields the
following constraint
\begin{align}
\label{eq:constraint-ua}
\langle\Delta^2\tilde u_{ij},\tilde u_{ij}\rangle-2|\nabla\tilde u_{ij}|^4=0.
\end{align}
Applying the identities provided by Lemma \ref{lem:der-u}
we find that this constraint is equivalent to
\begin{align*}
8m\frac{\sin^2\alpha}{r^4}\big(2-m+m\sin^2\alpha)=0
\end{align*}
yielding
\begin{align*}
2-m+m\sin^2\alpha=0.
\end{align*}
Now, assuming that the above constraint is satisfied, 
for the first \(m^2-1\) components of \(\tilde u_{ij}\), we are 
left with
\begin{align*}
\Delta^2\tilde u_{ij}
+2\operatorname{div}\big(|\nabla\tilde u_{ij}|^2\nabla\tilde u_{ij}\big)=0.
\end{align*}
Finally, a direct calculation using the identities given in \eqref{identities-u1} shows that
the first \(m^2-1\) components of \(\tilde u_{ij}\) satisfy 
\begin{align*}
\operatorname{div}\big(|\nabla\tilde u_{ij}|^2\nabla\tilde u_{ij}\big)
=&-\big(\frac{4m^2}{r^4}\sin^2\alpha\big)\tilde u_{ij},\\
\Delta^2\tilde u_{ij}=&\frac{8m(m-2)}{r^4}
\tilde u_{ij}.
\end{align*}
In conclusion, we get that for the first \(m^2-1\) components
\begin{align*}
\Delta^2\tilde u_{ij}
+2\operatorname{div}\big(|\nabla\tilde u_{ij}|^2\nabla\tilde u_{ij}\big)
=\frac{8m}{r^4}\big(m-2-m\sin^2\alpha\big)\tilde u_{ij}
\end{align*}
leading to the same condition \(2-m+m\sin^2\alpha=0\) such
that the existence part of the proof is now complete.

In order to show that \(\tilde u_{ij}\colon B^m\to \s^{m^2-1}\) is a weak biharmonic map we have to check if \(\tilde u_{ij}\in W^{2,2}(B^m,\s^{m^2-1})\).
We find that
\begin{align*}
\int_{B^m}|\nabla\tilde u_{ij}|^2\dv    
=&2(m-2)\vol(\s^{m-1})\int_0^1r^{m-3}dr,\\
\int_{B^m}|\Delta\tilde u_{ij}|^2\dv    
=&4m(m-2)\vol(\s^{m-1})\int_0^1r^{m-5}dr.
\end{align*}
We realize that we need to require \(m\geq 5\) in order for the second integral to be finite,
hence \(\tilde u_{ij}\) belongs to \(W^{2,2}(B^m,\s^{m^2-1})\) whenever \(m\geq 5\).
\end{proof}

As a second step we will investigate the stability of the 
proper biharmonic map provided by Theorem \ref{thm:biharmonic-u}.

\begin{proof}[Proof of Theorem \ref{thm:biharmonic-u-stability}]
We follow the ideas used in the proof of Theorem 1.2 of \cite{MR4076824}.
In order to prove the claim we will explicitly construct a variational vector field
for which the second variation of the bienergy, evaluated on this particular vector field, will be negative.

To this end, we consider a variation of \(\tilde u_{ij}\) defined as follows
\begin{align*}
\tilde u_{ij,s}=\big(\sin(\alpha+sV(r))u_{ij},\cos(\alpha+sV(r))\big),
\end{align*}
where \(s\in\R\) and \(V(r)\colon [0,1]\to [0,1]\) is a smooth function 
that will be specified during the proof. 
We need to require that \(V(1)=V'(1)=0\) in order to keep track of the boundary
conditions.

Recall the following identities
\begin{align*}
|\nabla u_{ij}|^2=\frac{2m}{r^2}, \qquad
|\Delta u_{ij}|^2=\frac{4m^2}{r^4}.
\end{align*}
Then, a direct calculation shows that
\begin{align*}
|\nabla \tilde u_{ij,s}|^4=s^4V'^4(r)+4m^2\sin^4(\alpha+sV(r))\frac{1}{r^4}
+4ms^2\sin^2(\alpha+sV(r))\frac{V'^2(r)}{r^2}.
\end{align*}
Moreover, we find
\begin{align*}
\Delta \tilde u_{ij,s} 
=&\bigg(-\sin(\alpha+sV(r))s^2 V'^2(s) u_{ij}
+\cos(\alpha+sV(r))s\Delta V(r) u_{ij} \\
&+\sin(\alpha+sV(r))\Delta u_{ij}, \\
&-\cos(\alpha+sV(r))s^2V'^2(r)-\sin(\alpha+sV(r))s \Delta V(r)
\bigg),
\end{align*}
where we used that \(x_l\nabla_l u_{ij}=0\)
which follows from a direct calculation.

Thus, we can infer
\begin{align*}
|\Delta \tilde u_{ij,s}|^2=&s^4V'^4(r)+s^2|\Delta V(r)|^2
+\sin^2(\alpha+sV(r))4m\big(\frac{m}{r^4}+s^2\frac{V'^2(r)}{r^2}\big)\\
&-\sin(\alpha+sV(r))\cos(\alpha+sV(r))
4s\frac{m}{r^2}\Delta V(r)
\end{align*}
such that we find
\begin{align*}
|\Delta \tilde u_{ij,s}|^2-|\nabla \tilde u_{ij,s}|^4=&
s^2|\Delta V(r)|^2
-\sin(\alpha+sV(r))\cos(\alpha+sV(r))4s\frac{m}{r^2}\Delta V(r)\\
&+\sin^2(\alpha+sV(r))\cos^2(\alpha+sV(r))4\frac{m^2}{r^4}\\
=&\big(s\Delta V(r)
-\sin(2\alpha+2sV(r))\frac{m}{r^2}
\big)^2.
\end{align*}

Hence, we obtain
\begin{align*}
\frac{d^2}{ds^2}\big|_{s=0}E_2(\tilde u_{ij,s})
=&\int_{B^m}\big(\Delta V(r)-2m\cos(2\alpha)\frac{V(r)}{r^2}\big)^2 \dv
-4m^2\sin^2(2\alpha)\int_{B^m}\frac{V^2(r)}{r^4}\dv.
\end{align*}
Also, by a direct calculation we get
\begin{align*}
\Delta V(r)=V''(r)+(m-1)\frac{V'(r)}{r}  
\end{align*}
and together with the condition for being proper biharmonic, which is
\(\sin^2\alpha=1-\frac{2}{m}\), we arrive at
\begin{align}
\label{eq:sv-u-ode-a}
\frac{d^2}{ds^2}\big|_{s=0}E_2(\tilde u_{ij,s})
=&\int_{B^m}\big(
V''(r)+(m-1)\frac{V'(r)}{r}+2(m-4)\frac{V(r)}{r^2}\big)^2 \dv \\
\nonumber&+32(2-m)\int_{B^m}\frac{V^2(r)}{r^4}\dv.
\end{align}

Now, we set \(V(r):=(1-r^2)^p,p>2\), which satisfies \(V(1)=V'(1)=0\) as requested.
Using this choice in \eqref{eq:sv-u-ode-a} we arrive, after some direct calculations, at the following expression
\begin{align*}
\frac{1}{\vol(\s^{m-1})}\frac{d^2}{ds^2}\big|_{s=0}E_2(\tilde u_{ij,s})
=&16p^2(p-1)^2\int_0^1(1-r^2)^{2p-4}r^{m+3}dr \\
&-16mp^2(p-1)\int_0^1(1-r^2)^{2p-3}r^{m+1}dr \\
&+4p\big(m^2p+4(p-1)(m-4)\big)\int_0^1(1-r^2)^{2p-2}r^{m-1}dr \\
&-8mp(m-4)\int_0^1(1-r^2)^{2p-1}r^{m-3}dr \\
&+4\big(m^2-16m+32\big)\int_0^1(1-r^2)^{2p}r^{m-5}dr.
\end{align*}

For \(a,b>0\) we recall the following integral formula
\begin{align*}
\int_0^1(1-x^2)^ax^bdx=\frac{\Gamma(a+1)\Gamma(\frac{b+1}{2})}{2\Gamma(a+\frac{b}{2}+\frac{3}{2})},   \end{align*}
where \(\Gamma(x)\) represents the Gamma function, which leads us to
\begin{align*}
\frac{\Gamma(2p-1+\frac{m}{2})}{8\vol(\s^{m-1})}\frac{d^2}{ds^2}\big|_{s=0}E_2(\tilde u_{ij,s})
=&4p^2(p-1)^2\Gamma(\frac{m}{2}+2)\Gamma(2p-3) \\
&-4mp^2(p-1)\Gamma(\frac{m}{2}+1)\Gamma(2p-2) \\
&+p\big(m^2p+4(p-1)(m-4)\big)\Gamma(\frac{m}{2})\Gamma(2p-1) \\
&-2mp(m-4)\Gamma(\frac{m}{2}-1)\Gamma(2p) \\
&+\big(m^2-16m+32\big)\Gamma(\frac{m}{2}-2)\Gamma(2p+1).
\end{align*}
Now, our strategy is as follows: For \(m\in\N\) given we need to find a \(p\)
such that the right hand side of the above equation is negative.
Using a computer algebra system one can directly check that for \(5\leq m\leq 12\) one can take 
\(p=m\) to obtain the claim completing the proof.
\end{proof}

Now, we turn to the proof of Theorem \ref{thm:biharmonic-v}.
First, we establish the following technical Lemma.

\begin{Lem}
\label{lem:der-v}
Let \(v_{ijk}\colon B^m\to\s^{m^3-1}\) be the map defined in \eqref{eqn:v}.
Then, the following identities hold
\begin{align}
|\nabla v_{ijk}|^2&=3\frac{m+1}{r^2},\\
\nonumber\Delta v_{ijk}&=\frac{3(m+1)}{\sqrt{m(m+2)}}
\big(-\frac{1}{r}\nabla_ku_{ij}+m\frac{x_k}{r^3}u_{ij}\big)\\
\nonumber&=-3\frac{m+1}{r^2}v_{ijk},\\
\nonumber\Delta^2 v_{ijk}&=\frac{15}{r^4}(m+1)(m-1)v_{ijk},
\end{align}
where \(u_{ij}\) is the map defined in \eqref{eqn:u}.
\end{Lem}
\begin{proof}
First of all, we note that we can express \(v_{ijk}\)
in terms of \(u_{ijk}\) as follows
\begin{align*}
v_{ijk}=\frac{1}{\sqrt{m(m+2)}}\big(r\nabla_ku_{ij}-m\frac{x_k}{r}u_{ij})
\end{align*}
such that we can make use of the identities obtained in Lemma \ref{lem:der-u}.

Again, a direct calculation shows that
\begin{align*}
\nabla_av_{ijk}=&\frac{1}{\sqrt{m(m+2)}}
\big(
\frac{x_a}{r}\nabla_ku_{ij}+r\nabla_a\nabla_ku_{ij}
-m\frac{\delta_{ak}}{r}u_{ij}
+m\frac{x_kx_a}{r^3}u_{ij}
-m\frac{x_k}{r}\nabla_au_{ij}
\big)
\end{align*}
and we can deduce
\begin{align*}
|\nabla v_{ijk}|^2=
\frac{1}{m(m+2)}\big(
&(1+m)^2|\nabla u_{ij}|^2+r^2|\nabla^2u_{ij}|^2+\frac{m^2(m-1)}{r^2} \\
&-\frac{4m}{r^2}|\underbrace{x\cdot\nabla u_{ij}}_{=0}|^2
+(1-m)\underbrace{x\cdot\nabla|\nabla u|^2}_{=-4m/r^2}
\big).
\end{align*}
Using the identity
\begin{align*}
|\nabla^2u_{ij}|^2=\frac{1}{2}\Delta|\nabla u_{ij}|^2
-\langle\nabla\Delta u_{ij},\nabla u_{ij}\rangle    
\end{align*}
we find \(|\nabla^2u_{ij}|^2=\frac{2m}{r^4}(4+m)\)
and by combining the equations we obtain the first claim of the Lemma.

Differentiating again we find
\begin{align*}
\nabla_b\nabla_av_{ijk}=\frac{1}{\sqrt{m(m+2)}}
\big(
&\frac{\delta_{ab}}{r}\nabla_ku_{ij}-\frac{x_ax_b}{r^3}\nabla_ku_{ij}
+\frac{x_a}{r}\nabla_b\nabla_ku_{ij}\\
&+\frac{x_b}{r}\nabla_a\nabla_ku_{ij}+r\nabla_b\nabla_a\nabla_ku_{ij} \\
&+\frac{mx_b\delta_{ak}}{r^3}u_{ij}-\frac{m\delta_{ak}}{r}\nabla_bu_{ij}\\
&+m\frac{\delta_{kb}x_a}{r^3}u_{ij}+m\frac{x_k\delta_{ab}}{r^3}u_{ij}
-3m\frac{x_kx_ax_b}{r^5}u_{ij}+m\frac{x_ax_k}{r^3}\nabla_bu_{ij} \\
&-m\frac{\delta_{kb}}{r}\nabla_au_{ij}
+m\frac{x_kx_b}{r^3}\nabla_au_{ij}
-m\frac{x_k}{r}\nabla_b\nabla_au_{ij}
\big).
\end{align*}
Taking the trace then gives
\begin{align*}
\Delta v_{ijk}=\frac{1}{\sqrt{m(m+2)}}
\big(  
&-\frac{m+1}{r}\nabla_ku_{ij}
+2\frac{x_b}{r}\nabla_b\nabla_ku_{ij}
+m(m-1)\frac{x_k}{r^3}u_{ij}\\
&+r\nabla_k\Delta u_{ij}
-m\frac{x_k}{r}\Delta u_{ij}
\big).
\end{align*}
Employing the identity \(x_b\nabla_k\nabla_bu_{ij}=-\nabla_ku_{ij}\)
and using the formula for \(\Delta u_{ij}\) derived in Lemma \ref{lem:der-u} we obtain the expression for \(\Delta v_{ijk}\).

Finally, using the product rule for the Laplace operator, we obtain
\begin{align*}
\Delta^2 v_{ijk}=
6\frac{(m+1)(m-4)}{r^4}v_{ijk}
+12\frac{m+1}{r^4}x_a\nabla_a v_{ijk}
-3\frac{m+1}{r^2}\Delta v_{ijk}.
\end{align*}
A direct calculation shows that \(x_a\nabla_a v_{ijk}=0\), together 
with the formula for \(\Delta v_{ijk}\) derived previously the proof is complete.
\end{proof}

\begin{proof}[Proof of Theorem \ref{thm:biharmonic-v}]
We use similar arguments as in the proof of Theorem \ref{thm:biharmonic-u}.
Again, we employ the following version of the equation
for biharmonic maps to spheres
\begin{align*}
\Delta^2u+2\operatorname{div}\big(|\nabla u|^2\nabla u\big)    
-\big(\langle\Delta^2u,u\rangle-2|\nabla u|^4\big)u=0.
\end{align*}
Since the \(m^3\)-th component of the map \(\tilde v_{ijk}\) is constant, 
we find the following constraint
\begin{align*}
\langle\Delta^2\tilde v_{ijk},\tilde v_{ijk}\rangle-2|\nabla\tilde v_{ijk}|^4=0.
\end{align*}
Inserting the identities obtained in Lemma \ref{lem:der-v}
we find that this constraint is equivalent to
\begin{align*}
3(m+1)\frac{\sin^2\beta}{r^4}\big(5m-5-(6m+6)\sin^2\beta\big)=0
\end{align*}
yielding
\begin{align*}
5m-5-(6m+6)\sin^2\beta=0.
\end{align*}
Now, assuming that the above constraint is satisfied, 
for the first \(m^3-1\) components of \(\tilde v_{ijk}\), we get
\begin{align*}
\Delta^2\tilde v_{ijk}
+2\operatorname{div}\big(|\nabla\tilde v_{ijk}|^2\nabla\tilde v_{ijk}\big)=0.
\end{align*}
Finally, a direct calculation using \eqref{identities-u1} shows that
the first \(m^3-1\) components of \(\tilde v_{ijk}\) satisfy the following identities
\begin{align*}
\operatorname{div}\big(|\nabla\tilde v_{ijk}|^2\nabla\tilde v_{ijk}\big)
=|\nabla\tilde v_{ijk}|^2\Delta \tilde v_{ijk}
=\big(-9\frac{(m+1)^2}{r^4}\sin^2\beta\big)\tilde v_{ijk}.
\end{align*}
In conclusion, we get that
\begin{align*}
\Delta^2\tilde v_{ijk}
+2\operatorname{div}\big(|\nabla\tilde v_{ijk}|^2\nabla\tilde v_{ijk}\big)
=3\frac{m+1}{r^4}\big(5m-5-(6m+6)\sin^2\beta\big)\tilde v_{ijk}
\end{align*}
leading to the same condition \(5m-5-(6m+6)\sin^2\beta=0\) showing that \(\tilde v_{ijk}\)
solves the biharmonic map equation except at the origin.
To complete the proof we need to show that
\(\tilde v_{ijk}\colon B^m\to \s^{m^3-1}\) is a weak biharmonic map, that is we have to check when \(\tilde v_{ijk}\in W^{2,2}(B^m,\s^{m^3-1})\).
We find that
\begin{align*}
\int_{B^m}|\nabla\tilde v_{ijk}|^2\dv    
=&5\frac{m-1}{2}\vol(\s^{m-1})\int_0^1r^{m-3}dr,\\
\int_{B^m}|\Delta\tilde v_{ijk}|^2\dv    
=&15\frac{m^2-1}{2}\vol(\s^{m-1})\int_0^1r^{m-5}dr.
\end{align*}
We realize that, again, we need to require \(m\geq 5\) in order for the second integral to be finite completing the proof.
\end{proof}

\begin{proof}[Proof of Theorem \ref{thm:biharmonic-v-stability}]
We use a similar strategy as in the proof of Theorem \ref{thm:biharmonic-u-stability}.
Again, consider a variation of \(\tilde v_{ijk}\) defined as follows
\begin{align*}
\tilde v_{ijk,s}=\big(\sin(\beta+sV(r))v_{ijk},\cos(\beta+sV(r))\big),
\end{align*}
where \(s\in\R\) and \(V(r)\colon\R\to\R\) is a function that will be fixed later.

Then, a direct calculation shows that
\begin{align*}
|\nabla \tilde v_{ijk,s}|^4=s^4V'^4(r)+\sin^4(\beta+sV(r))|\nabla v_{ijk}|^4
+2s^2\sin^2(\beta+sV(r))V'^2(r)|\nabla v_{ijk}|^2
\end{align*}
and also
\begin{align*}
\Delta \tilde v_{ijk,s} 
=&\bigg(-\sin(\beta+sV(r))s^2 V'^2(s) v_{ijk}
+\cos(\beta+sV(r))s\Delta V(r) v_{ijk} \\
&+\sin(\beta+sV(r))\Delta v_{ijk}, \\
&-\cos(\beta+sV(r))s^2V'^2(r)-\sin(\beta+sV(r))s \Delta V(r)
\bigg),
\end{align*}
where we used that \(x_l\nabla_l v_{ijk}=0\)
which follows from a direct calculation.

Hence, we can deduce that
\begin{align*}
|\Delta \tilde v_{ijk,s}|^2-|\nabla \tilde v_{ijk,s}|^4
=&\bigg(s\Delta V(r)
-\frac{1}{2}\sin(2\beta+2sV(r))|\nabla v_{ijk}|^2
\bigg)^2,
\end{align*}
where we used that \(|\Delta v_{ijk}|^2=|\nabla v_{ijk}|^4\) which holds 
as \(v_{ijk}\) is a harmonic map with values in the sphere.
Hence, we obtain
\begin{align*}
\frac{d^2}{ds^2}&\big|_{s=0}E_2(\tilde v_{ijk,s})
=\int_{B^m}\big(
\Delta V(r)-\cos(2\beta)V(r)|\nabla v_{ijk}|^2\big)^2 \dv
-\sin^2(2\beta)\int_{B^m}|\nabla v_{ijk}|^4V^2(r)\dv.
\end{align*}

Recall the following identity 
\begin{align*}
|\nabla v_{ijk}|^2=3\frac{m+1}{r^2},
\end{align*}
which was derived in Lemma \ref{lem:der-v},
and, also
\begin{align*}
\Delta V(r)=V''(r)+(m-1)\frac{V'(r)}{r}.    
\end{align*}
Note that the condition for \(\tilde v_{ijk}\) being proper biharmonic is
\(\sin^2\beta=\frac{5m-5}{6m+6}\), leading to
\begin{align*}
\cos 2\beta=\frac{2}{3}\frac{4-m}{m+1},\qquad    
\sin^2 2\beta=\frac{5}{9}\frac{(m-1)(m+11)}{(m+1)^2}.
\end{align*}
Using this data in the formula for the second variation we get
\begin{align}
\label{eq:sv-v-ode-a}
\frac{d^2}{ds^2}\big|_{s=0}E_2(\tilde v_{ijk,s})
=&\int_{B^m}\big(
V''(r)+(m-1)\frac{V'(r)}{r}++2(m-4)\frac{V(r)}{r^2}\big)^2 \dv \\
\nonumber&-5(m-1)(m+11)\int_{B^m}\frac{V^2(r)}{r^4}\dv.
\end{align}
Again, we now employ \(V(r)=(1-r^2)^p,p>2\), which satisfies \(V(1)=V'(1)=0\) as requested.
Using this choice in \eqref{eq:sv-v-ode-a} we get 
\begin{align*}
\frac{1}{\vol(\s^{m-1})}\frac{d^2}{ds^2}\big|_{s=0}E_2(\tilde v_{ijk,s})
=&16p^2(p-1)^2\int_0^1(1-r^2)^{2p-4}r^{m+3}dr \\
&-16mp^2(p-1)\int_0^1(1-r^2)^{2p-3}r^{m+1}dr \\
&+4p\big(m^2p+4(p-1)(m-4)\big)\int_0^1(1-r^2)^{2p-2}r^{m-1}dr \\
&-8mp(m-4)\int_0^1(1-r^2)^{2p-1}r^{m-3}dr \\
&-(m^2+82m-119)\int_0^1(1-r^2)^{2p}r^{m-5}dr.
\end{align*}

Carrying out the integrals as in the proof of Theorem \ref{thm:biharmonic-u-stability}
we find 
\begin{align*}
\frac{\Gamma(2p-1+\frac{m}{2})}{2\vol(\s^{m-1})}\frac{d^2}{ds^2}\big|_{s=0}E_2(\tilde v_{ijk,s})
=&16p^2(p-1)^2\Gamma(\frac{m}{2}+2)\Gamma(2p-3) \\
&-16mp^2(p-1)\Gamma(\frac{m}{2}+1)\Gamma(2p-2) \\
&+4p\big(m^2p+4(p-1)(m-4)\big)\Gamma(\frac{m}{2})\Gamma(2p-1) \\
&-8mp(m-4)\Gamma(\frac{m}{2}-1)\Gamma(2p) \\
&-(m^2+82m-119)\Gamma(\frac{m}{2}-2)\Gamma(2p+1).
\end{align*}
Now, for \(5\leq m\leq 18\) we can choose \(p=m\) to show the instability
of the biharmonic map \(\tilde v_{ijk}\) using a computer algebra system.
\end{proof}

\textbf{Data Availability Statement}: Data sharing not applicable to this article as no datasets were generated or
analysed during the current study.

\bibliographystyle{plain}

\begin{thebibliography}{10}

\bibitem{MR4599666}
Rare\c{s} Ambrosie, Cezar Oniciuc, and Ye-Lin Ou.
\newblock Biharmonic homogeneous polynomial maps between spheres.
\newblock {\em Results Math.}, 78(4):Paper No. 159, 40, 2023.

\bibitem{MR2044031}
Paul Baird and John~C. Wood.
\newblock {\em Harmonic morphisms between {R}iemannian manifolds}, volume~29 of
  {\em London Mathematical Society Monographs. New Series}.
\newblock The Clarendon Press, Oxford University Press, Oxford, 2003.

\bibitem{MR4106647}
V.~Branding, S.~Montaldo, C.~Oniciuc, and A.~Ratto.
\newblock Higher order energy functionals.
\newblock {\em Adv. Math.}, 370:107236, 60, 2020.

\bibitem{MR4040175}
Volker Branding and Yong Luo.
\newblock A nonexistence theorem for proper biharmonic maps into general
  {R}iemannian manifolds.
\newblock {\em J. Geom. Phys.}, 148:103557, 9, 2020.

\bibitem{MR3990379}
Volker Branding and Cezar Oniciuc.
\newblock Unique continuation theorems for biharmonic maps.
\newblock {\em Bull. Lond. Math. Soc.}, 51(4):603--621, 2019.

\bibitem{MR4477489}
Volker Branding and Anna Siffert.
\newblock On the equivariant stability of harmonic self-maps of cohomogeneity
  one manifolds.
\newblock {\em J. Math. Anal. Appl.}, 517(2):Paper No. 126635, 19, 2023.

\bibitem{MR1692148}
Sun-Yung~A. Chang, Lihe Wang, and Paul~C. Yang.
\newblock A regularity theory of biharmonic maps.
\newblock {\em Comm. Pure Appl. Math.}, 52(9):1113--1137, 1999.

\bibitem{MR164306}
James Eells, Jr. and J.~H. Sampson.
\newblock Harmonic mappings of {R}iemannian manifolds.
\newblock {\em Amer. J. Math.}, 86:109--160, 1964.

\bibitem{MR4076824}
A.~Fardoun, S.~Montaldo, and A.~Ratto.
\newblock Weakly biharmonic maps from the ball to the sphere.
\newblock {\em Geom. Dedicata}, 205:167--175, 2020.

\bibitem{MR4124860}
Hideaki Fujioka.
\newblock An example of harmonic map into the spheres with the singularity of
  order 4.
\newblock {\em J. Geom. Phys.}, 156:103810, 10, 2020.

\bibitem{MR581329}
M.~Giaquinta and J.~Ne\v{c}as.
\newblock On the regularity of weak solutions to nonlinear elliptic systems of
  partial differential equations.
\newblock {\em J. Reine Angew. Math.}, 316:140--159, 1980.

\bibitem{MR1401417}
Wenge Hao, Salvatore Leonardi, and Jind\v{r}ich Ne\v{c}as.
\newblock An example of irregular solution to a nonlinear {E}uler-{L}agrange
  elliptic system with real analytic coefficients.
\newblock {\em Ann. Scuola Norm. Sup. Pisa Cl. Sci. (4)}, 23(1):57--67, 1996.

\bibitem{MR214004}
Philip Hartman.
\newblock On homotopic harmonic maps.
\newblock {\em Canadian J. Math.}, 19:673--687, 1967.

\bibitem{MR1078114}
Fr\'{e}d\'{e}ric H\'{e}lein.
\newblock R\'{e}gularit\'{e} des applications faiblement harmoniques entre une
  surface et une sph\`ere.
\newblock {\em C. R. Acad. Sci. Paris S\'{e}r. I Math.}, 311(9):519--524, 1990.

\bibitem{MR714341}
Stefan Hildebrandt.
\newblock Nonlinear elliptic systems and harmonic mappings.
\newblock In {\em Proceedings of the 1980 {B}eijing {S}ymposium on
  {D}ifferential {G}eometry and {D}ifferential {E}quations, {V}ol. 1, 2, 3
  ({B}eijing, 1980)}, pages 481--615. Sci. Press Beijing, Beijing, 1982.

\bibitem{MR886529}
Guo~Ying Jiang.
\newblock {$2$}-harmonic maps and their first and second variational formulas.
\newblock {\em Chinese Ann. Math. Ser. A}, 7(4):389--402, 1986.
\newblock An English summary appears in Chinese Ann. Math. Ser. B {{\bf{7}}}
  (1986), no. 4, 523.

\bibitem{MR2375314}
Yin~Bon Ku.
\newblock Interior and boundary regularity of intrinsic biharmonic maps to
  spheres.
\newblock {\em Pacific J. Math.}, 234(1):43--67, 2008.

\bibitem{MR4227593}
Paul Laurain and Longzhi Lin.
\newblock Energy convexity of intrinsic bi-harmonic maps and applications {I}:
  {S}pherical target.
\newblock {\em J. Reine Angew. Math.}, 772:53--81, 2021.

\bibitem{MR916327}
Fang-Hua Lin.
\newblock A remark on the map {$x/|x|$}.
\newblock {\em C. R. Acad. Sci. Paris S\'{e}r. I Math.}, 305(12):529--531,
  1987.

\bibitem{MR2135286}
E.~Loubeau and C.~Oniciuc.
\newblock The index of biharmonic maps in spheres.
\newblock {\em Compos. Math.}, 141(3):729--745, 2005.

\bibitem{MR4371934}
Masashi Misawa and Nobumitsu Nakauchi.
\newblock Two examples of harmonic maps into spheres.
\newblock {\em Adv. Geom.}, 22(1):23--31, 2022.

\bibitem{MR3357596}
S.~Montaldo, C.~Oniciuc, and A.~Ratto.
\newblock Rotationally symmetric biharmonic maps between models.
\newblock {\em J. Math. Anal. Appl.}, 431(1):494--508, 2015.

\bibitem{MR4593065}
Nobumitsu Nakauchi.
\newblock A family of examples of harmonic maps into the sphere with one point
  singularity.
\newblock {\em Ex. Countex.}, 3:Paper No. 100107, 4, 2023.

\bibitem{MR4686854}
Nobumitsu Nakauchi.
\newblock Instability of a family of examples of harmonic maps.
\newblock {\em Ann. Global Anal. Geom.}, 65(1):10, 2024.

\bibitem{MR566246}
J.~Ne\v{c}as, O.~John, and J.~Star\'{a}.
\newblock Counterexample to the regularity of weak solution of elliptic
  systems.
\newblock {\em Comment. Math. Univ. Carolin.}, 21(1):145--154, 1980.

\bibitem{MR4265170}
Ye-Lin Ou and Bang-Yen Chen.
\newblock {\em Biharmonic submanifolds and biharmonic maps in {R}iemannian
  geometry}.
\newblock World Scientific Publishing Co. Pte. Ltd., Hackensack, NJ, [2020]
  \copyright 2020.

\bibitem{MR2285745}
Tristan Rivi\`ere.
\newblock Conservation laws for conformally invariant variational problems.
\newblock {\em Invent. Math.}, 168(1):1--22, 2007.

\bibitem{MR2383929}
Tristan Rivi\`ere and Michael Struwe.
\newblock Partial regularity for harmonic maps and related problems.
\newblock {\em Comm. Pure Appl. Math.}, 61(4):451--463, 2008.

\bibitem{MR375386}
R.~T. Smith.
\newblock The second variation formula for harmonic mappings.
\newblock {\em Proc. Amer. Math. Soc.}, 47:229--236, 1975.

\bibitem{MR2413109}
Michael Struwe.
\newblock Partial regularity for biharmonic maps, revisited.
\newblock {\em Calc. Var. Partial Differential Equations}, 33(2):249--262,
  2008.

\bibitem{MR2094320}
Changyou Wang.
\newblock Remarks on biharmonic maps into spheres.
\newblock {\em Calc. Var. Partial Differential Equations}, 21(3):221--242,
  2004.

\end{thebibliography}

\end{document}